\documentclass[12pt, oneside]{amsart}
\usepackage{geometry}
\numberwithin{equation}{section}                		
\usepackage{graphicx}				
\usepackage{amsmath,amssymb,amsthm,amscd}
\usepackage{ascmac}
\usepackage{indentfirst}
\usepackage{amsthm}
\usepackage{longtable}
\usepackage{color}
\usepackage{cases}
\usepackage[all]{xy}
\usepackage[OT2,T1]{fontenc}
\DeclareSymbolFont{cyrletters}{OT2}{wncyr}{m}{n}
\DeclareMathSymbol{\Sha}{\mathalpha}{cyrletters}{"58}
\DeclareMathOperator{\Gal}{\mathrm{Gal}}
\newtheorem{thm}{\textbf{Theorem}}[section]
\newtheorem*{thm*}{Theorem}
\newtheorem{mainthm}{Theorem}

\newtheorem{cor}[thm]{\textbf{Corollary}}
\newtheorem*{cor*}{\textbf{Corollary}}

\newtheorem*{conj*}{\textbf{Conjecture}}
\newtheorem{defn}[thm]{\textbf{Definition}}
\newtheorem{prop}[thm]{\textbf{Proposition}}
\newtheorem*{prop*}{\textbf{Proposition}}
\newtheorem{lem}[thm]{\textbf{Lemma}}

\newtheorem{que}[thm]{\textbf{Question}}
\theoremstyle{definition}
\newtheorem{example}[thm]{\textbf{Example}}
\newtheorem{rem}[thm]{Remark}

\renewcommand{\proofname}{\textit{Proof.}}

\newcommand{\Q}[1][]{\mathbb{Q}}
\newcommand{\Z}[1][]{\mathbb{Z}}
\newcommand{\C}[1][]{\mathbb{C}}
\newcommand{\F}[1][]{\mathbb{F}}
\newcommand{\Sel}[1][]{\mathrm{Sel}}
\newcommand{\Ker}[1][]{\mathrm{Ker}}
\newcommand{\Res}[1][]{\mathrm{Res}}
\newcommand{\Hom}[1][]{\mathrm{Hom}}
\newcommand{\loc}[1][]{\mathrm{Loc}}
\newcommand{\Fr}[1][]{\mathrm{Frob}}
\newcommand{\cl}[1][]{\mathrm{Cl}}
\newcommand{\ur}[1][]{\mathrm{ur}}
\newcommand{\A}[1][]{\mathcal{A}}
\newcommand{\mcO}[1][]{\mathcal{O}}
\newcommand{\tens}[1]{%
\mathbin{\mathop{\otimes}\limits_{#1}}%
}
\usepackage{hyperref}
\title[]{Ideal class groups of division fields of elliptic curves and everywhere unramified rational points}

\author{Naoto Dainobu}
\address{Department of Mathematics \\
3-14-1 Hiyoshi, Kohoku-ku, Yokohama-shi, Kanagawa 223-8522 Japan}
\email{vicarious@keio.jp}
					
\begin{document}
\maketitle
\begin{abstract}
Let $E$ be an elliptic curve over $\Q$, $p$ an odd prime number and $n$ a positive integer. In this article, we investigate the ideal class group $\cl(\Q(E[p^n]))$ of the $p^n$-division field $\Q(E[p^n])$ of $E$. We introduce a certain subgroup $E(\Q)_{\ur, p^n}$ of $E(\Q)$ and study the $p$-adic valuation of the class number $\#\cl(\Q(E[p^n]))$. 

In addition, when $n=1$, we further study $\cl(\Q(E[p]))$ as a $\Gal(\Q(E[p])/\Q)$-module. More precisely, we study the semi-simplification $(\cl(\Q(E[p]))\otimes \Z_p)^{\mathrm{ss}}$ of $\cl(\Q(E[p]))\otimes \Z_p$ as a $\Z_p[\Gal(\Q(E[p])/\Q)]$-module. We obtain a lower bound of the multiplicity of the $E[p]$-component in the semi-simplification when $E[p]$ is an irreducible $\Gal(\Q(E[p])/\Q)$-module.

\end{abstract}
\section{introduction}
\subsection{backgrounds}\label{back}
Ideal class groups of number fields have drawn many number theorists' interest for a long time. Among them, Kummer accomplished a monumental work on ideal class groups of cyclotomic fields in the 19th century toward Fermat's Last Theorem. Kummer studied the ideal class group $\cl(\Q(\mu_p))$ of the $p$-th cyclotomic field $\Q(\mu_p)$, where $p$ is an odd prime number and $\mu_p$ the group of $p$-th roots of unity. He examined the $p$-divisibility of the class number $\#\cl(\Q(\mu_p))$, or the non-vanishing of the $p$-primary part $A_p:=\cl(\Q(\mu_p)) \otimes \Z_p$ of $\cl(\Q(\mu_p))$. After that, Herbrand and Ribet refined Kummer's work in \cite{H} and \cite{R}. In their study, $A_p$ was treated with action of the Galois group $\Gal(\Q(\mu_p)/\Q)$. They decomposed $A_p$ with respect to this action and studied the non-vanishing of some components in the decomposition. 

In this article, as a natural analogue of cyclotomic fields, we consider division fields of elliptic curves and study their ideal class groups. Here we explain slight more details. Let $p$ be an odd prime number, $n$ a positive integer and $E$ an elliptic curve over $\Q$. We write 
$\Q(E[p^n])$ for the $p^n$-division field of $E$ which is generated over $\Q$ by the coordinates of all $p^n$-torsion points of $E$. We consider the following question on the ideal class group $\cl(\Q(E[p^n]))$ of $\Q(E[p^n])$. 
\begin{que}\label{q1}
When we have the $p$-divisibility $p | \#\cl(\Q(E[p^n]))$? Moreover, how large is the $p$-adic valuation of $\#\cl(\Q(E[p^n]))$?  
\end{que}

Furthermore, when $n=1$, we study the $p$-primary part $A(E)_p:=\cl(\Q(E[p]))\otimes \Z_p$ of $\cl(\Q(E[p]))$ with the action of the Galois group $\Gal(\Q(E[p])/\Q)$ as in Herbrand and Ribet's work. We write $A(E)_p^{\mathrm{ss}}$ for the semi-simplification of $A(E)_p$ as a $\Z_p[\Gal(\Q(E[p])/\Q)]$-module. Roughly speaking, $A(E)_p^{\mathrm{ss}}$ is a kind of a decomposition of $A(E)_p$ into irreducible $\Gal(\Q(E[p])/\Q)$-modules. Namely,  
\begin{align}\label{Pdecomp}
A(E)_p^{\mathrm{ss}} = \bigoplus_{M} M^{\oplus r_M},
\end{align}
where $M$ in the above direct product runs through all irreducible $\F_p[\Gal(\Q(E[p])/\Q)]$-modules and the non-negative integer $r_M$ denotes the multiplicity of the $M$-component. In this setting, we also consider the following question.

\begin{que}\label{q2}
For each irreducible component $M$ in (\ref{Pdecomp}), how large is the multiplicity $r_M$?
\end{que}

Let us introduce some known results on the above questions here. Concerning Question \ref{q1}, several lower bounds of the $p$-adic valuation of $\#\cl(\Q(E[p^n]))$ were given in \cite{SY1}, \cite{SY2}, \cite{Hi} and \cite{Oh} roughly as 
\begin{align}\label{pre}
v_p(\#\cl(\Q(E[p^n]))) \geqslant 2n (r(E)-1)-(\text{local contribution}\geqslant 0).
\end{align}
Here $v_p$ denotes the additive $p$-adic valuation normalized by $v_p(p)=1$ and $r(E)$ denotes the Mordell-Weil rank of $E/\Q$. By the previous results (\ref{pre}), we can obtain non-trivial lower bounds of $v_p(\#\cl(\Q(E[p^n])))$ when $r(E) \geqslant 2$ in general. While, the so-called ``minimalist conjecture'' implies that asymptotically zero percent of all elliptic curves over $\Q$ have rank at least $2$. Concerning this conjecture, see \cite[Section 1]{BMSW} and \cite[Section 1.1]{BHKSSW}. Thus it seems important to consider Question \ref{q1} when $r(E) \leqslant 1$.

A recent work by Prasad and Shekhar \cite{PS} greatly motivated us to consider Question \ref{q2}. Using the $p$-Selmer group $\Sel(\Q, E[p])$ of $E$, they studied the non-vanishing of the $E[p]$-component in $A(E)_p^{\mathrm{ss}}$ when $E[p]$ is an irreducible $\Gal(\Q(E[p])/\Q)$-module. In \cite[Theorem 3.1]{PS}, they proved an implication
\begin{align}\label{Pimp}
\dim_{\F_p}(\Sel(\Q, E[p]))\geqslant 2 \Rightarrow r_{E[p]} \neq 0
\end{align}
under certain conditions. When $p\geqslant 11$, the condition $\dim_{\F_p}(\Sel(\Q, E[p]))\geqslant 2$ is equivalent to the condition $r(E)\geqslant 2$ or $\Sha(E/\Q)[p]\neq 0$, where $\Sha(E/\Q)[p]$ denotes the $p$-torsion part of the Tate-Shafarevich group $\Sha(E/\Q)$ of $E$. The interest in \cite{PS} is mainly in the case where  $\Sha(E/\Q)[p]\neq 0$. While, it seems rare that we have $\Sha(E/\Q)[p]\neq 0$ for an odd prime number $p$. As we noted above, the ``minimalist conjecture'' implies that we rarely have $r(E) \geqslant 2$. Therefore, to consider Question \ref{q2}, the author thinks it is important to improve (\ref{Pimp}) so that we can treat the case where $\dim_{\F_p}(\Sel(\Q, E[p])) \leqslant 1$.

\subsection{Main results}\label{intro}
 
In our main results, we give some partial answers to the two questions in the previous subsection. In the following, we write 
\[
K_{E, p^n}:=\Q(E[p^n])
\]
 for the $p^n$-division field of $E$ and $h_{E, p^n}:=\#\cl(K_{E, p^n})$ for its class number.

  We first introduce a certain subgroup $E(\Q)_{\ur, p^n}$ of $E(\Q)$ which we call {\it everywhere unramifed rational points}.
\begin{defn}\label{unrampt}
We define a subgroup $E(\Q)_{\ur, p^n}$ of $E(\Q)$ as 
\[
E(\Q)_{\ur, p^n}:=\Ker\left(E(\Q) \rightarrow \prod_{\ell : \text{{\rm prime}}} \frac{E(\Q_{\ell}^{\ur})}{p^nE(\Q_{\ell}^{\ur})}\right),
\]
where $\ell$ in the above product runs over all prime numbers and $\Q_{\ell}^{\ur}$ denotes the maximal unramified extension of $\Q_{\ell}$. We further define a non-negative integer
\[
r_{\ur, p^n}(E):=l_{\Z_p}(E(\Q)_{\ur, p^n}/p^nE(\Q)),
\]
where for a $\Z_p$-module $M$, $l_{\Z_p}(M)$ denotes its length as a $\Z_p$-module.
\end{defn}
In short, $E(\Q)_{\ur, p^n}$ is the group of points in $E(\Q)$ which are $p^n$-divisible in $E(\Q^{\ur}_{\ell})$ for every prime number $\ell$. The non-negative integer $r_{\ur, p^n}(E)$ measures the difference between the group $p^nE(\Q)$ of globally $p^n$-divisible points and the group $E(\Q)_{\ur, p^n}$ of everywhere locally $p^n$-divisible points.

The group of everywhere unramifed rational points $E(\Q)_{\ur, p^n}$ is related to the ideal class group $\cl(K_{E, p^n})$ as follows. 
\begin{mainthm}[Theorem \ref{main1}]\label{A}
Suppose that the condition
\begin{itemize}
\item[(Inj)] $H^1(K_{E, p^n}/\Q, E[p^n])=0$
\end{itemize}
holds. Then
\[
l_{\Z_p}\left(\Hom_{\Gal(K_{E, p^n}/\Q)}(\cl(K_{E, p^n}), E[p^n])\right)\geqslant r_{\ur, p^n}(E).
\]
Moreover, if we further assume the irreducibility of $E[p]$ as a $\Gal(K_{E, p}/\Q)$-module, then we obtain 
\[
v_p(h_{E, p^n})\geqslant 2r_{\ur, p^n}(E).
\] 
\end{mainthm}
\begin{rem}\label{13}
\begin{enumerate}
\item[(1)]In \cite{LW}, for $p\geqslant 3$, Lawson and Wuthrich gave a necessary and sufficient condition for the above condition (Inj). By their result, we can see that (Inj) always holds for any $n$ when $p\geqslant13$.
\item[(2)] Suppose that $\mathrm{rank}_{\Z}(E(\Q)_{\ur, p^n})=r(E)$ holds in addition to (Inj) and the irreducibility of $E[p]$, where $\mathrm{rank}_{\Z}(E(\Q)_{\ur, p^n})$ denotes the $\Z$-rank of the torsion-free part of $E(\Q)_{\ur, p^n}$. Then we obtain an inequality
\[
v_p(h_{E, p^n}) \geqslant 2r_{\ur, p^n}(E) \geqslant 2nr(E)
\]
by the latter statement in Theorem \ref{A}. This gives a partial improvement of the lower bounds (\ref{pre}) in previous results, and we can obtain a non-trivial lower bound of $v_p(h_{E, p^n})$ when $\mathrm{rank}_{\Z}(E(\Q)_{\ur, p^n})=r(E) \geqslant 1$. Later, we introduce such an example in Example \ref{ex} in Section 5. In addition, we can obtain a non-trivial lower bound of $v_p(h_{E, p^n})$ even when $r(E) = 0$, using Theorem \ref{A} in some cases. Such a case is treated in Example \ref{Hi}.
\end{enumerate}
\end{rem}

Therefore, when $E[p]$ is irreducible, the value $r_{\ur, p^n}(E)$ gives us a partial answer to Question \ref{q1}. When $n=1$, the value $r_{\ur, p}(E)$ further gives a partial answer to Question \ref{q2} for $M=E[p]$.

\begin{cor*}[Corollary \ref{cor1}]\label{corA}
Suppose the irreducibility of $E[p]$ as a $\Gal(K_{E, p}/\Q)$-module. Then there exists a $\Gal(K_{E, p}/\Q)$-equivariant surjective homomorphism 
\[
(\cl(K_{E, p})\tens{\Z}{\F_p})^{\mathrm{ss}}\twoheadrightarrow E[p]^{\oplus r_{\ur, p}(E)}.
\]
Here $(\cl(K_{E, p})\otimes\F_p)^{\mathrm{ss}}$ denotes the semi-simplification of $\cl(K_{E, p})\otimes\F_p$ as a $\Gal(K_{E, p}/\Q)$-module. In particular, we obtain
\[
r_{E[p]} \geqslant r_{\ur, p}(E),
\]
where $r_{E[p]}$ is the non-negative integer defined in (\ref{Pdecomp}).
\end{cor*}
\begin{rem}
\begin{enumerate}
\item[(1)] The latter part of the above corollary follows by the existence of  an injection 
\[
(\cl(K_{E,p})\otimes\F_p)^{\mathrm{ss}} \hookrightarrow A(E)_p^{\mathrm{ss}}=(\cl(K_{E,p})\otimes\Z_p)^{\mathrm{ss}}.
\]  
\item[(2)]
This corollary gives a partial improvement of (\ref{Pimp}) by Prasad and Shekhar. In fact, we can show that 
\[
\dim_{\F_p}(\Sel(\Q, E[p]))\geqslant 2 \Rightarrow r_{\ur, p}(E) \geqslant 1
\]
when $\Sha(E/\Q)[p]=0$ and $E(\Q_p)[p]=0$ hold. Using this corollary, we can obtain the non-vanishing $r_{E[p]} \neq 0$ for $E$ which satisfies $\dim_{\F_p}(\Sel(\Q, E[p]))=1$ in some cases. We treat such a case in Example \ref{ex}. 
\end{enumerate}
\end{rem}

Thus for our Questions \ref{q1} and \ref{q2}, it is important to study the value $r_{\ur, p^n}(E)$ which is determined by the group $E(\Q)_{\ur, p^n}$. Concerning $E(\Q)_{\ur, p^n}$, we have the following.
\begin{prop*}[Proposition \ref{mainprop}]
Suppose that the condition
\begin{itemize}
\item[(Tam)] For any prime number $\ell$ $(\neq p)$, we have $v_p(c(E/\Q_{\ell}))=0$. Here $v_p$ denotes the normalized $p$-adic valuation and $c(E/\Q_{\ell})$ the Tamagawa number of $E/\Q_{\ell}$.
\end{itemize}
holds. Then we have 
\[
E(\Q)_{\ur, p^n}=\Ker\left(E(\Q) \rightarrow E(\Q_p^{\ur})/p^nE(\Q_p^{\ur})\right).
\]
\end{prop*}
In other words, if the condition (Tam) in the proposition hold, then $E(\Q)_{\ur, p^n}$ is determined by the local information on $E(\Q)$ only at the fixed prime number $p$. Although this is an easy consequence of \cite[Lemma 2.5 (3)]{PS}, we deduce this proposition by a general property of Tamagawa factors for $p$-adic representations in Section 3. This proposition gives an useful sufficient condition on a rational point $P \in E(\Q)$ of infinite order for $P \in E(\Q)_{\ur, p^n}$.

\begin{mainthm}[Theorem \ref{main2}]\label{B}
Assume that the condition $\mathrm{(Tam)}$ holds and $E$ is defined by a Weierstrass equation which is minimal at $p$. If there exists a rational point $P:=(X, Y) \in E(\Q)$ such that 
\begin{center}
$v_p(X)<0$\ \ and \ \ $v_p\left(X/Y\right) \geqslant n+1$,
\end{center}
then $P \in E(\Q)_{\ur, p^n}$. 
\end{mainthm}
\begin{rem}
Due to the Nagell-Lutz theorem \cite[Chapter VIII, Corollary 7.2]{Sil}, the rational point $P$ in Theorem \ref{B} has infinite order. Thus Theorem \ref{B} can be applicable to $E$ which satisfies $r(E)\geqslant 1$. 
\end{rem}

We prove these main results in Section 4 after introducing basic properties of $E(\Q)_{\ur, p^n}$ and Selmer groups in Sections 2 and 3. Theorem \ref{A} is a consequence of the properties of $E(\Q)_{\ur, p^n}$ and Ohshita's result in \cite[Lemma 2.10]{Oh}. By Proposition \ref{mainprop}, the proof of Theorem \ref{B} is reduced to the investigation of the local Mordell-Weil group $E(\Q_p)$ with the formal group logarithm attached to $E/\Q_p$. After the proofs of the main results, we give some numerical examples in Section 5.

\subsection*{Acknowledgement}
The author would like to thank his supervisor Masato Kurihara for his continued support and careful reading of a draft of this paper. The author would like to thank Yoshinosuke Hirakawa and Hideki Matsumura for helpful discussions and generous support for this work. Thanks are also due to Takenori Kataoka for valuable comments on Lemma \ref{ktok}. This research was supported by JSPS KAKENHI Grant Number 21J13502.

\section{prerilinaries}\label{prelim}
\setcounter{subsection}{-1}
\subsection{Notation}
Let $F$ be a number field or a local field. We write $G_F$ for its absolute Galois group $\Gal(\overline{F}/F)$, where $\bar{F}$ is a fixed algebraic closure of $F$. We put $F^{\ur}$ as the maximal unramified abelian extension of $F$ in $\overline{F}$. For  a number field $F$ and a place $v$ of $F$, $F_v$ denotes the completion of $F$ at $v$, and $I_v$ denotes the inertia subgroup of $v$ in the decomposition group $G_{F_v}$. 

Let $L/F$ be a Galois extension of number fields or local fields and $N$ a $\Gal(L/F)$-module. We abbreviate the Galois cohomology group $H^i(\Gal(L/F), N)$ by $H^i(L/F, N)$. In particular, when $L=\bar{F}$, we further abbreviate $H^i(G_F, N)$ by $H^i(F, N)$. For a number field $F$ and a place $v$ of $F$, we define the unramified subgroup of the Galois cohomology group $H^i(F_v, N)$ as
\[
H_{\ur}^i(F_v, N):=
\Ker\left(H^i(F_v, N)\xrightarrow{\Res} H^i(F^{\ur}_v, N) \right).
\]
 Here $\Res$ denotes the restriction of cohomology classes in $H^i(F_v, N)$ to the inertia subgroup $I_v$ of $v$.

\subsection{Everywhere unramified rational points}
In the following, we fix an odd prime number $p$, a positive integer $n$, an elliptic curve $E$ over $\Q$. We abbreviate $K_{E, p^n}=\Q(E[p^n])$ by $K_n$ if no confusion occurs.

First, we relate $E(\Q)_{\ur, p^n}$ in Definition \ref{unrampt} to the ideal class group $\cl(K_n)$ of $K_n$. We consider the $p^n$-th Kummer map $\kappa_n$ and the restriction $\Res_{K_n/\Q}$ of cohomology classes
\[
E(\Q)/p^nE(\Q)\xrightarrow{\kappa_n} H^1(\Q, E[p^n]) \xrightarrow{\Res_{K_n/\Q}} \Hom_{\Gal(K_n/\Q)}(G_{K_n}, E[p^n]).
\]
Here we note that 
\[
\mathrm{Im} (\Res_{K_n/\Q}) \subset H^1(K_n, E[p^n])^{\Gal(K_n/\Q)}=\Hom_{\Gal(K_n/\Q)}(G_{K_n}, E[p^n])
\]
 since $G_{K_n}$ acts trivially on $E[p^n]$.

\begin{prop}\label{factor}
We have an inclusion
\[
\Res_{K_n/\Q}\circ  \kappa_n \left(E(\Q)_{\ur, p^n}/p^nE(\Q) \right) \subset \Hom_{\Gal(K_n/\Q)}(\cl(K_n), E[p^n]).
\]
\end{prop}
\proofname\ \ For every place $v$ in $\Q$ and every place $w$ in $K_n$ such that $w\mid v$, we have a commutative diagram  
\[
  \xymatrix{
      E(\Q)/p^nE(\Q) \ar[r]^{\kappa_n}\ar[d]^{incl} & H^1(\Q, E[p^n])\ar[d]^{\loc^{\ur}_v}\ar[r]^-{\hspace{1mm}\Res_{K_n/\Q}} & \hspace{2mm}\Hom_{\Gal(K_n/\Q)}(G_{K_n}, E[p^n])\ar[d]\\
  E(\Q^{\ur}_v)/p^nE(\Q^{\ur}_v)\ar[r]^{\kappa_{n, v}} &  H^1(\Q^{\ur}_v,  E[p^n])\ar[r]^-{\Res_{w/v}} & \Hom(I_{w}, E[p^n]).
  }
\]
Here $\Res_{w/v}$ is the restriction of cohomology classes to $I_{w}$, the right vertical arrow is the restriction of homomorphisms from $G_{K_n}$ to $I_{w}$, and $\kappa_{n, v}$ is the local Kummer map at $v$. 

Suppose $v$ is a finite place of $\Q$. By the definition of the $E(\Q)_{\ur, p^n}$, the image of $E(\Q)_{\ur, p^n}/p^nE(\Q)$ under the left vertical arrow is zero, which implies $\Res_{K_n/\Q}\circ \kappa_n(x)$ is trivial when restricted to $I_{w}$ for any element $x \in E(\Q)_{\ur, p^n}/p^nE(\Q)$. If $v$ is the infinite place of $\Q$, then $\Hom(I_{w}, E[p^n])=0$ since we assume that $p$ is an odd prime number. Thus the image of $E(\Q)_{\ur, p^n}/p^nE(\Q)$ under $\Res_{K_n/\Q}\circ \kappa_n$ are homomorphisms which factor through $\Gal(H/K_n)$, where $H$ denotes the Hilbert class field of $K_n$. Due to class field theory, we have an isomorphism $\cl(K_n)\simeq \Gal(H/K_n)$ as $\Gal(K_n/\Q)$-modules, which implies the proposition. \hfill$\square$

\subsection{Selmer group}\label{sel}
In this subsection, we introduce some basic properties of Selmer groups in the Bloch-Kato style. The definitions and results provided here are more general than those in other parts of this article. Here the author would like to thank Ryotaro Sakamoto who taught him some basic properties of Tamagawa factors.

Let $p$ be a prime number and $V$ a $p$-adic representation of $G_{\Q}$ over $\Q_p$. We take a $G_{\Q}$-stable $\Z_p$-lattice $T$ of $V$. We define discrete $G_{\Q}$-modules 
\[
A:=V/T\ \ \text{and} \ \ M_n:=T/p^nT\cong A[p^n]\ \ (n \in \Z_{>0}).
\]
 First, we recall the definition of the Bloch-Kato Selmer group. For more details, see \cite[Section 5]{BK}.

\begin{defn}
For every place $v$ of $\Q$, Bloch and Kato's local condition $H_f^1(\Q_v, V)$ in $H^1(\Q_v, V)$ is defined as  
\[
   \begin{cases}
    H_f^1(\Q_v, V):=  H_{\ur}^1(\Q_v, V)=\Ker\left(H^1(\Q_v, V)\rightarrow H^1(\Q^{\ur}_v, V)\right)   
    & (v \neq p), \\
    H_f^1(\Q_p, V):= \Ker\left(H^1(\Q_p, V)\rightarrow H^1(\Q_p, V \otimes\mathbf{B}_{\mathrm{crys}})\right) & (v=p).
  \end{cases}
\]
Here $\mathbf{B}_{\mathrm{crys}}$ denotes Fontaine's crystalline period ring, which is defined in \cite[Section2]{Fo} and \cite[Section 1]{BK}, for example. 
\end{defn}

We have exact sequences 
\begin{align}\label{TVA}
0\rightarrow T\xrightarrow{\iota} V \xrightarrow{\pi} A \rightarrow 0,
\end{align}
\begin{align}\label{pAA}
0 \rightarrow M_n \xrightarrow{i_n} A \xrightarrow{\times p^n} A \rightarrow 0.
\end{align}  
By $\iota$, $\pi$ and $i_n$, we have the induced homomorphisms 
\[
\iota : H^1(\Q_v, T) \rightarrow H^1(\Q_v, V),\ \ \pi : H^1(\Q_v, V) \rightarrow H^1(\Q_v, A)
\]
\[
i_n : H^1(\Q_v, M_n) \rightarrow H^1(\Q_v, A).
\]
\begin{defn}\label{loccondiTAM}
We define local conditions at $v$ with coefficients in $T$, $A$ and $M_n$ as 
\[
H^1_f(\Q_v, T):=\iota^{-1}(H^1_f(\Q_v, V)),\ \ H^1_f(\Q_v, A):=\pi(H^1_f(\Q_v, V))
\]
\[
H^1_f(\Q_v, M_n):=i_n^{-1}(H^1_f(\Q_v, A))
\]
respectively.
\end{defn}
 \begin{defn}\label{Sel}
Let $W$ be one of $V$, $T$, $A$ and $M_n$. The Bloch-Kato Selmer group of $W$ is defined as
\[
H^1_f(\Q, W) := \Ker\left(H^1(\Q, W) \xrightarrow{\prod \mathrm{Loc}_v} \prod_{v : \text{place}} \frac{H^1(\Q_v, W)}{H^1_f(\Q_v, W)} \right).
\]
Here $v$ in the above product runs over all places in $\Q$ and $\mathrm{Loc}_v$ denotes the restriction map $H^1(\Q, W)\rightarrow H^1(\Q_v, W)$ at a place $v$.
\end{defn}

\begin{rem}\label{exell}
For an elliptic curve $E$ over $\Q$, let $V:=V_pE$ be the rational $p$-adic Tate module of $E$ and $T:=T_p E$ be the integral $p$-adic Tate module of $E$. We have $A \cong E[p^{\infty}]$ and $M_n\cong E[p^n]$. The Bloch-Kato Selmer group of $W=V, T, A, M_n$ is the same as the classical Selmer group $\Sel(\Q, W)$ of $E$. Here we check this fact briefly. It suffices to show that $H^1_f(\Q_v, V)$ coincides with the classical local condition $H^1_{\mathrm{cl}}(\Q_v, V)$ in $H^1(\Q_v, V)$ for every place $v$ of $\Q$. This is because the classical local conditions with coefficients in $T, A, M_n$ are determined by $H^1_{\mathrm{cl}}(\Q_v, V)$  in exactly the same way as Definition \ref{loccondiTAM}. The local condition $H^1_{\mathrm{cl}}(\Q_v, V)$ is defined as the image of the Kummer map 
\[
\kappa : E(\Q_v)^{\wedge}\tens{\Z_p}{\Q_p} \rightarrow H^1_f(\Q_v, V). \vspace{-3mm}
\]
Here, $E(\Q_v)^{\wedge}$ is the $p$-adic completion of $E(\Q_v)$. The equality $H^1_f(\Q_v, V) = \mathrm{Im}{\kappa}$ follows from the first commutative diagram in \cite[Section 3, Example 3.10.1]{BK}.
\end{rem}

For a prime number $\ell$ $(\neq p)$, we note that $H^1_f(\Q_{\ell}, M_n)$ is not necessarily equal to the unramified subgroup $H^1_{\ur}(\Q_{\ell}, M_n)$. This difference can be measured by the ($p$-part of the) Tamagawa factor of $A$ at $\ell$.
\begin{defn}\label{Tamdef}
For a prime number $\ell$ $(\neq p)$, we put 
\[
\A_{\ell}:=A^{I_{\ell}}/(A^{I_{\ell}})^{\mathrm{div}}.
\]
Here $(A^{I_{\ell}})^{\mathrm{div}}$ denotes the maximal divisible subgroup of $A^{I_{\ell}}$. We define the $p$-part of the Tamagawa factor $c_{A, \ell}$ of $A$ at $\ell$ as 
\[
c_{A, \ell}:=\#\A_{\ell}/(\Fr_{\ell}-1)=\#\A_{\ell}^{\Fr_{\ell}=1},
\]
where $\Fr_{\ell}$ denotes a Frobenius element in $\Gal(\Q_{\ell}^{\ur}/\Q_{\ell})$.
\end{defn}
\begin{rem}
This definition of Tamagawa factors for $p$-adic representations coincides with the one defined in \cite[Section 4]{DSW}.
\end{rem}

\begin{prop}\label{Tam}
If $c_{A, \ell}=1$, then $H^1_f(\Q_{\ell}, M_n)=H^1_{\ur}(\Q_{\ell}, M_n)$.
\end{prop}
\proofname\ We first consider the local condition with the coefficient in $A$. \begin{align*}
H^1_{\ur}(\Q_{\ell}, A)/H^1_f(\Q_{\ell}, A) & \overset{\sim}{\longrightarrow}  \mathrm{coker}\left(H^1_f(\Q_{\ell}, V)=H^1_{\ur}(\Q_{\ell}, V) \overset{\pi}{\longrightarrow} H^1_{\ur}(\Q_{\ell}, A)\right)\\
&\overset{\sim}{\longrightarrow} \mathrm{coker}\left(V^{I_{\ell}}/(\mathrm{\Fr_{\ell}}-1)(V^{I_{\ell}}) \overset{\pi}{\longrightarrow} A^{I_{\ell}}/(\mathrm{\Fr_{\ell}}-1)(A^{I_{\ell}})\right)\\
&\overset{\sim}{\longrightarrow} \A_{\ell}/(\mathrm{\Fr_{\ell}}-1)\A_{\ell}.
\end{align*}
In the second isomorphism above, we use the following.
\begin{lem}{\cite[Appendix B, Lemma 2.8]{Ru}}\label{ruap}
Let $G$ be a pro-cyclic group and $\gamma$ a topological generator of $G$. Suppose that $W$ is a $\Z_p[G]$-module which is either a finitely generated $\Z_p$-module, or a finite dimensional $\Q_p$-vector space, or a discrete torsion $\Z_p$-module. Then we have an isomorphism 
\[
H^1(G, W) \xrightarrow[\mathrm{ev}]{\sim} W/(\gamma-1)W,
\]
where the isomorphism $\mathrm{ev}$ is induced by evaluating cocycles at $\gamma$.
\end{lem}
Thus if $c_{A, \ell}=1$, then we have $H^1_{\ur}(\Q_{\ell}, A)=H^1_f(\Q_{\ell}, A)$. On the other hand, by the definition of $H^1_f(\Q_{\ell}, M_n)$, we have  
\begin{align*}
H^1_f(\Q_{\ell}, M_n)&=\Ker(H^1(\Q_{\ell}, M_n) \rightarrow H^1(\Q_{\ell}, A) \rightarrow H^1(\Q^{\ur}_{\ell}, A)^{\Fr_{\ell}=1})\\
&= \Ker(H^1(\Q_{\ell}, M_n) \rightarrow H^1(\Q^{\ur}_{\ell}, M_n)^{\Fr_{\ell}=1} \overset{\phi}{\rightarrow} H^1(\Q^{\ur}_{\ell}, A)^{\Fr_{\ell}=1}).
\end{align*}
We can see that the above homomorphism $\phi$ is injective. In fact, by (\ref{pAA}), we have an exact sequence
\[
0\rightarrow (A^{I_{\ell}} \otimes \Z/p^n\Z)^{\Fr_{\ell}=1} \longrightarrow H^1(\Q^{\ur}_{\ell}, M_n)^{\Fr_{\ell}=1}\overset{\phi}{\longrightarrow} H^1(\Q^{\ur}_{\ell}, A)^{\Fr_{\ell}=1}.
\]
Since $A^{I_{\ell}} \otimes \Z/p^n\Z = \A_{\ell}\otimes \Z/p^n\Z$, 
\[
\Ker(\phi)=(A^{I_{\ell}} \otimes \Z/p^n\Z)^{\Fr_{\ell}=1}=(\A_{\ell} \otimes \Z/p^n\Z)^{\Fr_{\ell}=1}.
\]
We have an exact sequence 
\[
0 \rightarrow \A_{\ell} \xrightarrow[\times p]{} \A_{\ell} \rightarrow \A_{\ell}\otimes \Z/p\Z \rightarrow 0
\]
from which we obtain 
\[
0 \rightarrow \A_{\ell}^{\Fr_{\ell}=1} \xrightarrow[\times p]{} \A_{\ell}^{\Fr_{\ell}=1} \rightarrow (\A_{\ell}\otimes \Z/p\Z)^{\Fr_{\ell}=1} \rightarrow H^1(\Q^{\ur}_{\ell}/\Q_{\ell}, \A_{\ell})\rightarrow \cdots.
\]
By the assumption $c_{A, \ell}=1$, we have $\A_{\ell}^{\Fr_{\ell}=1}=0$. Using the isomorphism $\mathrm{ev}$ in Lemma \ref{ruap} and again the assumption $c_{A, \ell}=1$, we obtain
\[
H^1(\Q^{\ur}_{\ell}/\Q_{\ell}, \A_{\ell}) \xrightarrow[\mathrm{ev}]{\sim} \A_{\ell}/(\Fr_{\ell}-1)\A_{\ell}=0.
\]
This implies $\Ker(\phi)=(\A_{\ell} \otimes \Z/p^n\Z)^{\Fr_{\ell}=1}=0$. Thus the homomorphism $\phi$ is injective, and we have
\[
H^1_f(\Q_{\ell}, M_n)= \Ker(H^1(\Q_{\ell}, M_n) \rightarrow H^1(\Q^{\ur}_{\ell}, M_n)^{\Fr_{\ell}=1} )=H^1_{\ur}(\Q_{\ell}, M_n).
\]
\hfill$\square$

\section{Unramifiedness of $\Sel(\Q, E[p^n])$ outside $p$}
Let $p$ be an odd prime number, $n$ a positive integer and $E$ an elliptic curve over $\Q$ again. In this section, we study the unramifiedness of the Selmer group $\Sel(\Q, E[p^n])=H^1_f(\Q, E[p^n])$ outside $p$, applying the general results in Subsection \ref{sel} to the case where $V=V_pE$ and $T=T_pE$. 

Putting $A:=V/T$, we first compare the Tamagawa factor in Definition \ref{Tamdef} with the usual (geometric) Tamagawa number of $E$. Let $\ell$ be a prime number, $F$ an algebraic extension of $\Q_{\ell}$, and $\F$ the residue field of $F$. We put $\widetilde{E}(\F)$ as the group of $\F$-rational points on the reduced curve $\widetilde{E}$ of $E$ modulo $\ell$, and $\widetilde{E}_{\mathrm{ns}}(\F)$ as the group of points in $\widetilde{E}(\F)$ which are non-singular. By the reduction map modulo the maximal ideal of $F$
\[
\mathrm{red} : E(F) \to \widetilde{E}(\F),
\]
we define
\[
E_0(F):= \{ Q \in E(F) \mid \mathrm{red}(Q) \in \widetilde{E}_{\mathrm{ns}}(\F)\}
\]
and 
\[
E_1(F):=\Ker(E_0(F) \xrightarrow[]{\mathrm{red}} \widetilde{E}_{\mathrm{ns}}(\F))
\]
as usual.

Let $c(E/\Q_{\ell}):=E(\Q_{\ell})/E_0(\Q_{\ell})$ be the (geometric) Tamagawa number of $E$ at $\ell$. 

\begin{lem}\label{Tamlem}
When $\ell$ $\neq p$, we have
\[
E_0(\Q_{\ell}^{\ur})[p^{\infty}] \simeq \widetilde{E}^{\mathrm{ns}}(\bar{\F}_{\ell})[p^{\infty}].
\]
\end{lem}

\vspace{2mm}
\proofname\ \ We have an exact sequence
\begin{align}\label{lemuse}
0 \to E_1(\Q_{\ell}^{\ur}) \to  E_0(\Q_{\ell}^{\ur}) \xrightarrow[]{\mathrm{red}} \widetilde{E}^{\mathrm{ns}}(\bar{\F}_{\ell})\to 0. 
\end{align}
The group $E_1(\Q_{\ell}^{\ur})$ is isomorphic to $\Z_{\ell}^{\ur}$ by the theory of the formal group attached to $E$ over $\Z_{\ell}$, where $\Z_{\ell}^{\ur}$ is the ring of integers of $\Q_{\ell}^{\ur}$. Then $E_1(\Q_{\ell}^{\ur})$ is $p$-divisible since $\ell\neq p$, and we have an exact sequence 
\[
0 \to E_1(\Q_{\ell}^{\ur})[p^{\infty}] \to  E_0(\Q_{\ell}^{\ur})[p^{\infty}] \to \widetilde{E}^{\mathrm{ns}}(\bar{\F}_{\ell})[p^{\infty}] \to 0
\]
by taking the $p$-power torsion part of each term of (\ref{lemuse}). Since $E_1(\Q_{\ell}^{\ur}) \simeq \Z_{\ell}^{\ur}$, we have $E_1(\Q_{\ell}^{\ur})[p^{\infty}]=0$, which implies the lemma. \hfill\qed

\begin{prop}\label{TamTam}
For a prime number $\ell$ ($\neq p$), we have 
\[
c_{A, \ell} = p^{v_p(c(E/\Q_{\ell}))},
\]
where $v_p$ denotes the normalized $p$-adic valuation.
\end{prop}

\proofname\ \ When $E$ has good reduction at $\ell$, then $A^{I_{\ell}} = A$ by the criterion of Néron-Ogg-Shafarevich \cite[Chapter 7, Theorem 7.1]{Sil}. Then we have $c_{A, \ell}=1$ by definition. While, we also have $c(E/\Q_{\ell})=1$ by definition, and the desired equality holds.

We consider the case where $E$ has bad reduction at $\ell$ in the following. We have an exact sequence 
\begin{align}\label{tamex}
0 \to E_0(\Q_{\ell}^{\ur}) \to E(\Q_{\ell}^{\ur}) \to E(\Q_{\ell}^{\ur})/E_0(\Q_{\ell}^{\ur}) \to 0. 
\end{align}
The group $E_0(\Q_{\ell}^{\ur})$ is a $p$-divisible group. For this, it suffices to show that $\widetilde{E}^{\mathrm{ns}}(\bar{\F}_{\ell})$ is $p$-divisible since we have the exact sequence (\ref{lemuse}) and the $p$-divisibility of $E_1(\Q_{\ell}^{\ur})$. The $p$-divisibility of $\widetilde{E}^{\mathrm{ns}}(\bar{\F}_{\ell})$ follows from the isomorphism 
\begin{align}\label{redmordstr}
\widetilde{E}^{\mathrm{ns}}(\bar{\F}_{\ell})\simeq
\begin{cases}
\bar{\F}^{\times}_{\ell} & (\text{when $E$ has multiplicative reduction at $\ell$}) \\
\bar{\F}_{\ell} & (\text{when $E$ has additive reduction at $\ell$}).
\end{cases}
\end{align}
Thus $E_0(\Q_{\ell}^{\ur})$ is $p$-divisible, and  we have an exact sequence 
\begin{align}\label{tamex-p}
0 \to E_0(\Q_{\ell}^{\ur})[p^{\infty}] \to E(\Q_{\ell}^{\ur})[p^{\infty}]=A^{I_{\ell}} \to (E(\Q_{\ell}^{\ur})/E_0(\Q_{\ell}^{\ur}))[p^{\infty}] \to 0
\end{align}
by taking the $p$-power torsion part of each term of (\ref{tamex}). Taking the quotient of each term of (\ref{tamex-p}) by its maximal divisible subgroup, we obtain another exact sequence
\begin{align}\label{tamex-div}
 (E_0(\Q_{\ell}^{\ur})[p^{\infty}])/(E_0(\Q_{\ell}^{\ur})[p^{\infty}])^{\mathrm{div}} \to \mathcal{A}_{\ell} \to (E(\Q_{\ell}^{\ur})/E_0(\Q_{\ell}^{\ur}))[p^{\infty}] \to 0.
\end{align}
Here we note that $\#E(\Q_{\ell}^{\ur})/E_0(\Q_{\ell}^{\ur}) < \infty$ by \cite[Corollary 9.2 (b)]{Silad}. As we noted above, the group $E_0(\Q_{\ell}^{\ur})$ is $p$-divisible, and so is its $p$-power torsion subgroup $E_0(\Q_{\ell}^{\ur})[p^{\infty}]$. Then we have $\mathcal{A}_{\ell} \simeq (E(\Q_{\ell}^{\ur})/E_0(\Q_{\ell}^{\ur}))[p^{\infty}]$ and 
\begin{align}\label{finish}
\mathcal{A}^{\Fr_{\ell}=1}_{\ell} \simeq (E(\Q_{\ell}^{\ur})/E_0(\Q_{\ell}^{\ur}))^{\Fr_{\ell}=1}[p^{\infty}]
\end{align}
by (\ref{tamex-div}). By taking the $\Fr_{\ell}$-fixed part of each term of (\ref{tamex-p}), we obtain an exact sequence
\begin{align}\label{tamrel}
0 &\to E_0(\Q_{\ell})[p^{\infty}] \to E(\Q_{\ell})[p^{\infty}] \to (E(\Q_{\ell}^{\ur})/E_0(\Q_{\ell}^{\ur}))^{\Fr_{\ell}=1}[p^{\infty}] \\
 & \to H^1(\Q_{\ell}^{\ur}/\Q_{\ell}, E_0(\Q_{\ell}^{\ur})[p^{\infty}])\nonumber.
\end{align}
By Lemma \ref{Tamlem}, we have 
\begin{align}\label{defect}
H^1(\Q_{\ell}^{\ur}/\Q_{\ell}, E_0(\Q_{\ell}^{\ur})[p^{\infty}]) = H^1(\Q_{\ell}^{\ur}/\Q_{\ell}, \widetilde{E}^{\mathrm{ns}}(\bar{\F}_{\ell})[p^{\infty}]).
\end{align}
When $E$ has additive reduction at $\ell$, then $\widetilde{E}^{\mathrm{ns}}(\bar{\F}_{\ell})[p^{\infty}]=0$ by (\ref{redmordstr}), which implies $H^1(\Q_{\ell}^{\ur}/\Q_{\ell}, E_0(\Q_{\ell}^{\ur})[p^{\infty}])=0$. When $E$ has split multiplicative reduction at $\ell$, then $\widetilde{E}^{\mathrm{ns}}(\bar{\F}_{\ell})[p^{\infty}]\simeq \mu_{p^{\infty}}(\bar{\F}_{\ell})$ by (\ref{redmordstr}). Here $\mu_{p^{\infty}}(\bar{\F}_{\ell})$ denotes the group of all $p$-power roots of unity in $\bar{\F}_{\ell}$. We note that this isomorphism is $\Fr_{\ell}$-equivariant when $E$ has split multiplicative reduction. See \cite[Chapter 3, Proposition 2.5]{Sil} for details on this isomorphism, for example. Hence, 
\begin{align*}
H^1(\Q_{\ell}^{\ur}/\Q_{\ell}, \widetilde{E}^{\mathrm{ns}}(\bar{\F}_{\ell})[p^{\infty}]) \simeq H^1(\Q_{\ell}^{\ur}/\Q_{\ell}, \mu_{p^{\infty}}(\bar{\F}_{\ell})) &\simeq \mu_{p^{\infty}}(\bar{\F}_{\ell})/(\Fr_{\ell}-1)\mu_{p^{\infty}}(\bar{\F}_{\ell})\\
& \simeq (\Q_p/\Z_p)/(\ell-1)(\Q_p/\Z_p) = 0.
\end{align*}
Here we use Lemma \ref{ruap} in the second isomorphism above. When $E$ has non-split multiplicative reduction at $\ell$, $E$ acquires split multiplicative reduction over the quadratic unramified extension $F$ over $\Q_{\ell}$. We have the inflation-restriction exact sequence
\begin{align*}
 0&\to H^1(F/\Q_{\ell}, \widetilde{E}^{\mathrm{ns}}(\bar{\F}_{\ell})^{G_F}[p^{\infty}]) \to
 H^1(\Q_{\ell}^{\ur}/\Q_{\ell}, \widetilde{E}^{\mathrm{ns}}(\bar{\F}_{\ell})[p^{\infty}]) \\
 &\to H^1(\Q_{\ell}^{\ur}/F, \widetilde{E}^{\mathrm{ns}}(\bar{\F}_{\ell})[p^{\infty}])^{\Gal(F/\Q_{\ell})}
 \to H^2(F/\Q_{\ell}, \widetilde{E}^{\mathrm{ns}}(\bar{\F}_{\ell})^{G_F}[p^{\infty}]).
\end{align*}
Since $p$ is odd, $H^i(F/\Q_{\ell}, \widetilde{E}^{\mathrm{ns}}(\bar{\F}_{\ell})^{G_F}[p^{\infty}]) = 0$ for $i\in\Z_{>0}$. Thus we have 
\[
H^1(\Q_{\ell}^{\ur}/\Q_{\ell}, \widetilde{E}^{\mathrm{ns}}(\bar{\F}_{\ell})[p^{\infty}]) \simeq H^1(\Q_{\ell}^{\ur}/F, \widetilde{E}^{\mathrm{ns}}(\bar{\F}_{\ell})[p^{\infty}])^{\Gal(F/\Q_{\ell})}=0.
\]
Here we use the equality $H^1(\Q_{\ell}^{\ur}/F, \widetilde{E}^{\mathrm{ns}}(\bar{\F}_{\ell})[p^{\infty}])=0$ which follows by almost the same argument as in the split multiplicative case. By (\ref{defect}), we obtain the equality $H^1(\Q_{\ell}^{\ur}/\Q_{\ell}, E_0(\Q_{\ell}^{\ur})[p^{\infty}])=0$ and 
\[
(E(\Q_{\ell})/E_0(\Q_{\ell}))[p^{\infty}] \simeq (E(\Q_{\ell}^{\ur})/E_0(\Q_{\ell}^{\ur}))^{\Fr_{\ell}=1}[p^{\infty}]
\]
by (\ref{tamrel}). Then (\ref{finish}) implies
\[
c_{A, \ell} = \#\mathcal{A}^{\Fr_{\ell}=1}_{\ell}= \#(E(\Q_{\ell})/E_0(\Q_{\ell}))[p^{\infty}]= p^{v_p(c(E/\Q_{\ell}))}.
\]
\hfill\qed

\begin{cor}\label{unram}
For a prime number $\ell$ $(\neq p)$, suppose that a condition
\begin{itemize}
\item[(Tam)]\ \ \ $v_p(c(E/\Q_{\ell}))=0$
\end{itemize}
holds. Then the restriction
\[
\loc_{\ell}^{\ur} : \Sel(\Q, E[p^n]) \rightarrow H^1(\Q_{\ell}^{\ur}, E[p^n])
\]
is 0-map.
\end{cor}

\proofname\ \ It suffices to show that $H^1_f(\Q_{\ell}, E[p^n])=H^1_{\ur}(\Q_{\ell}, E[p^n])$. The condition (Tam) implies $c_{A, \ell} = 1$ by Proposition \ref{TamTam}. This further implies the equality $H^1_f(\Q_{\ell}, E[p^n])=H^1_{\ur}(\Q_{\ell}, E[p^n])$ by Proposition \ref{Tam}.\hfill\qed

\begin{rem}\label{tamj}
When $E$ does not have split multiplicative reduction at $\ell$, it is well-known that the (geometric) Tamagawa number $c(E/\Q_{\ell})$ is at most $4$ \cite[Section 7, Theorem 6.1]{Sil}. Then in this case, the condition (Tam) is automatically satisfied if $p\geqslant 5$. When $E$ has split multiplicative reduction at $\ell$, then we have $c(E/\Q_{\ell})=-v_{\ell}(j(E))$ \cite[Section 7, Theorem 6.1]{Sil}, where $j(E)$ denotes the $j$-invariant of $E$. Therefore, the condition (Tam) is equivalent that $v_p(v_{\ell}(j(E)))=0$ in this case.
\end{rem}

\begin{rem}
If the condition (Tam) in Corollary \ref{unram} holds, the local condition $H^1_f(\Q_{\ell}, E[p^n])$ often vanishes. Indeed, the local condition is isomorphic to $E(\Q_{\ell})/p^nE(\Q_{\ell})$ via Kummer map by definition. The condition (Tam) implies that 
\[
E(\Q_{\ell})/p^nE(\Q_{\ell}) =E_0(\Q_{\ell})/p^nE_0(\Q_{\ell}).
\]
 We have an exact sequence
\[
0 \to E_1(\Q_{\ell}) \to E_0(\Q_{\ell}) \to  \widetilde{E}^{\mathrm{ns}}(\F_{\ell}) \to 0.
\]
Since $ E_1(\Q_{\ell}) \simeq \Z_{\ell}$ is $p$-divisible, we obtain $E_0(\Q_{\ell})/p^nE_0(\Q_{\ell}) \simeq \widetilde{E}^{\mathrm{ns}}(\F_{\ell})/p^n\widetilde{E}^{\mathrm{ns}}(\F_{\ell})$ by the above exact sequence. Thus when $E$ has good reduction at $\ell$, we have $H^1_f(\Q_{\ell}, E[p^n])=0$ if and only if $p\nmid\#\widetilde{E}(\F_{\ell})$. When $E$ has additive reduction at $\ell$, $\widetilde{E}^{\mathrm{ns}}(\F_{\ell}) \simeq \F_{\ell}$. This is a $p$-divisible group, which implies $H^1_f(\Q_{\ell}, E[p^n])=0$. Suppose that $E$ has multiplicative reduction at $\ell$. Since $\F_{\ell}$ is not algebraically closed, we have an isomorphism
\begin{align*}\label{multstr}
\widetilde{E}^{\mathrm{ns}}(\F_{\ell})\simeq
\begin{cases}
\Ker(N : \F^{\times}_{\ell^2}\to\F^{\times}_{\ell}) & (\text{if $E$ has non-split multiplicative reduction at $\ell$}) \\
\F_{\ell}^{\times} & (\text{if $E$ has split multiplicative reduction at $\ell$}).
\end{cases}
\end{align*}
Here $\F_{\ell^2}$ denotes the quadratic extension of $\F_{\ell}$ and $N : \F^{\times}_{\ell^2}\to\F^{\times}_{\ell}$ denotes the norm map. For this isomorphism, see \cite[Chapter 3, Exercise 3.5]{Sil} or \cite[Proposition 1.8.1]{Co} for instance. We have 
$\Ker(N : \F^{\times}_{\ell^2}\to\F^{\times}_{\ell}) = \{x \in \F^{\times}_{\ell^2}\mid x^{\ell+1} =1\}$. Hence, when $E$ has non-split multiplicative reduction, we have $H^1_f(\Q_{\ell}, E[p^n])=0$ if and only if $\ell\not \equiv -1\pmod p$. When $E$ has split multiplicative reduction, we have $H^1_f(\Q_{\ell}, E[p^n])=0$ if and only if $\ell\not \equiv 1\pmod p$.  
\end{rem}

\section{Proof of the main results}
Now we prove our main results which we introduced in Subsection \ref{intro}. 

\subsection{Theorem \ref{A} and its corollary}\mbox{}
Here we state Theroem \ref{A} in Subsection \ref{intro} again. We abbreviate the $p^n$-division field $K_{E, p^n}:=\Q(E[p^n])$ of $E$ by $K_n$. 
\begin{thm}\label{main1}
Suppose that the condition
\begin{itemize}
\item[(Inj)] $H^1(K_n/\Q, E[p^n])=0$
\end{itemize}
holds. Then
\[
l_{\Z_p}\left(\Hom_{\Gal(K_n/\Q)}(\cl(K_n), E[p^n])\right)\geqslant r_{\ur, p^n}(E).
\]
Moreover, if we further assume the irreducibility of $E[p]$ as a $\Gal(K_1/\Q)$-module, then we obtain 
\[
v_p(h_{E, p^n})\geqslant 2r_{\ur, p^n}(E).
\] 
Here we recall that $h_{E, p^n}$ denotes the class number $\#\cl(K_n)$ of $K_n$.
\end{thm}

\proofname\ \  We have a homomorphism 
\[
E(\Q)_{\ur, p^n}/p^nE(\Q) \xrightarrow{\Res_{K_n/\Q}\circ \kappa_n}  \Hom_{\Gal(K_n/\Q)}(\cl(K_n), E[p^n]) 
\]
by Proposition \ref{factor}. We recall that $\Res_{K_n/\Q}$ denotes the restriction map 
\[
H^1(\Q, E[p^n]) \rightarrow H^1(K_n, E[p^n])^{\Gal(K_n/\Q)}, 
\]
and this is injective under the condition (Inj). This implies the first inequality in the proposition.

Next, suppose that $E[p]$ is irreducible as a $\Gal(K_1/\Q)$-module. By \cite[Lemma 5]{Ru2}, $E[p]$ is absolutely irreducible as a $\Gal(K_1/\Q)$-module since $p$ is odd. The following result by Ohshita is a key for us.
\begin{lem}{\cite[Lemma 2.10]{Oh}}\label{Ohshita}
Let $F$ be a number field and $T$ a free $\Z_p$-module of rank $d$ on which $G_F$ acts continuously and linearly. Let $W:=T\tens{\Z_p}\Q_p/\Z_p$ and 
\[
F_n := F(W[p^n])
\]
be a number field which corresponds to the kernel of the representation 
\[
G_F \rightarrow \mathrm{Aut}_{\Z/p^n\Z}(W[p^n]).
\]
For $n \in \Z_{\geqslant 0}$ and a $\Z_p$-submodule $M$ of $\Hom_{\Gal(F_n/F)}(G_{F_n}, W[p^n])$, we define $F_n(M)$ as the fixed field of $\bigcap_{h\in M} \Ker h$. Then we have 
\[
[F_n(M):F_n] = p^{dl_{\Z_p}(M)},
\]
where $l_{\Z_p}(M)$ denotes the length of $M$ as a $\Z_p$-module.
\end{lem}

We use Lemma \ref{Ohshita} for $F=\Q$, $T=T_p E$, $F_n=K_n$, and 
\[
M=\Hom_{\Gal(K_n/\Q)}(\cl(K_n), E[p^n])=\Hom_{\Gal(K_n/\Q)}(\Gal(K_n^{\ur}/K_n), E[p^n]),
\]
where $K_n^{\ur}$ is the maximal unramified abelian extension of $K_n$. The number field $K_n(M)$ is an intermediate field of the extension $K_n^{\ur}/K$ by definition. Then Lemma \ref{Ohshita} implies
\[
v_p(h_{E, p^n}) \geqslant v_p([K_n(M):K_n]) = v_p(2l_{\Z_p}(M)) =  v_p(2r_{\ur, p^n}(E)).
\]
\hfill$\square$

\begin{cor}\label{cor1}
Suppose $n=1$ and the irreducibility of $E[p]$ as a $\Gal(K_1/\Q)$-module. Then there exists a $\Gal(K_1/\Q)$-equivariant surjective homomorphism 
\[
(\cl(K_1)\tens{\Z}{\F_p})^{\mathrm{ss}}\twoheadrightarrow E[p]^{\oplus r_{\ur, p}(E)}.
\]
 In particular, we obtain
\[
r_{E[p]} \geqslant r_{\ur, p}(E),
\]
where $r_{E[p]}$ is the non-negative integer defined in (\ref{Pdecomp}).
\end{cor}

\proofname\ \ First, we note that the condition (Inj) in Theorem \ref{main1} holds by the irreducibility of $E[p]$ as a $\Gal(K_1/\Q)$-module. We use the following lemma.
\begin{lem}\label{ktok}
Let $G$ be a finite group and $M$, $N$ be $\F_p[G]$-modules such that they are finite dimensional over $\F_p$. Suppose $N$ is absolutely irreducible as an $\F_p[G]$-module. If $\dim_{\F_p}\left(\Hom_{G}(M, N)\right)\geqslant r$, then there exists a $G$-equivariant surjective homomorphism 
\[
M^{\mathrm{ss}} \rightarrow N^{\oplus r}.
\]
\end{lem}
We prove this lemma later. By Theorem \ref{main1}, we have 
\[
\dim_{\F_p}\left(\Hom_{\Gal(K_1/\Q)}(\cl(K_1)\otimes \F_p, E[p])\right)\geqslant r_{\ur, p}(E).
\]
 Applying Lemma \ref{ktok} to the case where $G=\Gal(K_1/\Q)$, $M=\cl(K_1)\otimes \F_p$, $N=E[p]$ and $r=r_{\ur, p}(E)$, we obtain Corollary \ref{cor1}. Here we use \cite[Lemma 5]{Ru2} again, which states that the irreducibility of $E[p]$ as an $\F_p[\Gal(K_1/\Q)]$-module implies the absolute irreducibility in our setting. \hfill$\square$
 
\vspace{4mm}
\noindent\textit{Proof of Lemma \ref{ktok}.}\ \ We prove the statement by induction on $r$. When $r=1$, we clearly have a surjective $G$-equivariant homomorphism $M \twoheadrightarrow N$ by the irreducibility of $N$, which implies the existence of desired homomorphism $M^{\mathrm{ss}}\twoheadrightarrow N$. Next, let $r$ be a positive integer. Then we have a surjective $G$-equivariant homomorphism $f : M \twoheadrightarrow N$ again. Putting $M_0:=\Ker f$, we obtain a decomposition
\begin{eqnarray}\label{decomp}
 M^{\mathrm{ss}}=M_0^{\mathrm{ss}}\oplus N. 
 \end{eqnarray}
 On the other hand, we have an exact sequence 
\[
0\rightarrow M_0 \rightarrow M \xrightarrow{f} N\rightarrow 0.
\]
Applying $\Hom_{G}(\ \cdot\ , N)$ to the above exact sequence, we obtain
\[
0 \rightarrow \Hom_{G}(N, N) \rightarrow \Hom_{G}(M, N) \rightarrow \Hom_{G}(M_0, N).
\]
By the absolute irreducibility of $N$ as an $\F_p[G]$-module, we have $\Hom_{G}(N, N)\cong\F_p$ and $\dim_{\F_p}(\Hom_{G}(M_0, N))\geqslant r-1$. Then there exists a surjective $G$-equivariant homomorphism 
\[
M_0^{\mathrm{ss}}\twoheadrightarrow N^{\oplus(r-1)}
\]
 by the induction hypothesis. Thus we have a desired surjective $G$-equivariant homomorphism $M^{\mathrm{ss}}\twoheadrightarrow N^{\oplus r}$ by the decomposition (\ref{decomp}). \hfill$\square$

\subsection{Discription of $E(\Q)_{\ur, p^n}$ and Theorem \ref{B}}
\begin{prop}\label{mainprop}
Suppose the condition 
\begin{itemize}
\item[(Tam)] For any prime number $\ell$ $(\neq p)$, we have $v_p(c(E/\Q_{\ell}))=0$. Here $v_p$ denotes the normalized $p$-adic valuation and $c(E/\Q_{\ell})$ the Tamagawa number of $E/\Q_{\ell}$.
\end{itemize}
holds. Then we have 
\[
E(\Q)_{\ur, p^n}=\Ker\left(E(\Q) \rightarrow E(\Q_p^{\ur})/p^nE(\Q_p^{\ur})\right).
\]
\end{prop}
\proofname\ \  Due to Corollary \ref{unram}, the restriction
\[
\loc_{\ell}^{\ur} : \Sel(\Q, E[p^n]) \rightarrow H^1(\Q_{\ell}^{\ur}, E[p^n])
\]
is 0-map for any prime number $\ell$ $(\neq p)$ under the condition (Tam). On the other hand, we have a commutative diagram
\[
  \xymatrix{
      0\ar[r]&E(\Q)/p^nE(\Q) \ar[r]^{\kappa_n}\ar[d]^{incl} & \Sel(\Q, E[p^n])\ar[d]^{\loc^{\ur}_{\ell}}\\
  0\ar[r]&E(\Q^{\ur}_{\ell})/p^nE(\Q^{\ur}_{\ell})\ar[r]^{\kappa_{n, \ell}} &  H^1(\Q^{\ur}_{\ell},  E[p^n]).
  }
\]
Then the left vertical arrow is the 0-map, which implies the proposition. \hfill$\square$

Here we state Theorem \ref{B} in Subsection \ref{intro} again.
\begin{thm}\label{main2}
Assume that the condition $({\rm Tam})$ holds and $E$ is defined by a Weierstrass equation which is minimal at $p$. If there exists a rational point $P:=(X, Y) \in E(\Q)$ such that 
\begin{center}
$v_p(X)<0$\ \ and \ \ $v_p\left(X/Y\right) \geqslant n+1$,
\end{center}
then $P \in E(\Q)_{\ur, p^n}$. 
\end{thm}

\proofname\ \ The group $E_1(\Q_p)$ which we defined in the beginning of Section 3 has an alternative description as 
\[
E_1(\Q_p)=\{ Q=(X, Y) \in E(\Q_p) \mid v_p(X)<0 \}.
\]
By this expression and the first assumption on $P$, we have $P \in E_1(\Q_p)$. We have an isomorphism 
\[
E_1(\Q_p) \cong \hat{E}(p\Z_p)\ \ \ (x, y) \mapsto -x/y.
\]
by \cite[Chapter 7, Proposition 2.2]{Sil}, where $\hat{E}$ denotes the formal group associated to $E$. Since $p$ is odd, the formal logarithm $\log_{\hat{E}}$ of $\hat{E}$ induces the isomorphism 
\[
\log_{\hat{E}} : \hat{E}(p\Z_p) \cong p\Z_p.
\]
 Thus we have the isomorphism 
 \[
 E_1(\Q_p) \xrightarrow{\sim} p\Z_p\ \ \ (x, y) \mapsto \log_{\hat{E}}(-x/y).
\]
 The assumption $v_p\left(X/Y\right) \geqslant n+1$ implies that $P \in p^n E_1(\Q_p) \subset p^nE(\Q_p)$. Then we obtain $P \in E(\Q)_{\ur, p^n}$ by  Proposition \ref{mainprop}. \hfill\qed

\section{Numerical examples}
Finally, we give some numerical examples of our main results. Here the author thanks Yoshinosuke Hirakawa again for informing him an interesting example (Example \ref{Hi}) which deals with an elliptic curve of Mordell-Weil rank 0. 
\begin{example}\label{Hi}
Let $p=5$ and $E$ be the elliptic curve defined by the equation 
\[
y^2+y=x^3-x^2
\]
 which has LMFDB \cite{LMFDB} label $11.a3$. We apply Theorem \ref{main1} to this example for $n=1$.

We know that $j(E)=-2^{12}\cdot 11^{-11}$ due to LMFDB, and (Tam) in Proposition \ref{mainprop} is satisfied by Remark \ref{tamj}. Using MAGMA \cite{MAGMA}, we can see that the quadratic twist of $E$ by $5$ does not have rational $5$-torsion points, which implies that (Inj) in Theorem \ref{main1} holds by \cite[Theorem1]{LW}. On the other hand, LMFDB tells us that 
\[
E(\Q)=E(\Q)[5]=\langle(0,0)\rangle
\]
 and that $E(\Q_{5^5})[5^2] \cong \Z/5^2\Z$, where $\Q_{5^5}$ denotes the unramified extension of $\Q_5$ of degree $5$. Thus we can see that the rational $5$-torsion point $(0, 0)$ is $5$-divisible in $E(\Q^{\ur}_5)$. Then $(0, 0) \in E(\Q)_{\ur, 5}$ by Proposition \ref{mainprop}, and we have $r_{\ur, 5}(E) \geqslant 1$. This implies that $r_{\ur, 5}(E) = 1$ since $r_{\ur, 5}(E) \leqslant \dim_{\F_5}(E(\Q)/5E(\Q))=1$. We obtain 
\[
\Hom_{\Gal(K_{E,5}/\Q)}(\cl(K_{E,5}), E[5])) \neq 0
\]
by Theorem \ref{main1}. In particular, we obtain $v_5(h_{E, 5}) > 0$. 
\end{example}

\begin{example}\label{ex}
Let $p=13$ and $E$ be the elliptic curve defined by the equation 
\[
y^2+y=x^3+x^2
\]
which has LMFDB label $43.a1$. We apply Theorem \ref{main1} to this example for $n=1, 2$.

We know that $j(E)=2^{12}\cdot 43^{-1}$ due to LMFDB, and and (Tam) in Proposition \ref{mainprop} is satisfied by Remark \ref{tamj}. By Remark \ref{13}, (Inj) in Theorem \ref{main1} is satisfied for any $n\in\Z_{>0}$ since $p=13$. LMFDB also tells us that $r(E)=1$ and that $(0, 0)$ is a generator of $E(\Q)$. Using MAGMA, we obtain 
\[
19\cdot (0, 0) = \left(\frac{-2^3\cdot 3^2\cdot 5 \cdot 11 \cdot 59 \cdot 61 \cdot 107}{{\bf 13}^6\cdot 37^2}, \frac{3^4\cdot 11^2 \cdot 17 \cdot 59^2 \cdot 173 \cdot 211}{{\bf 13}^9\cdot 37^3}\right).
\]
By Theorem \ref{main2}, we have $19\cdot (0, 0)\in E(\Q)_{\ur, p^n}$ for $n=1, 2$. Since $19\cdot (0, 0) \notin 13E(\Q)$, 
\[
r_{\ur, 13^n}(E) \geqslant n
\] 
for $n=1, 2$. LMFDB also says that $E[13]$ is irreducible as a $\Gal(K_{E, 13}/\Q)$-module. Thus we obtain
\[
 v_{13}(h_{E, 13^n}) \geqslant 2n 
 \]
for $n=1, 2$ by Theorem \ref{main1}.

Furthermore, when $n=1$, there exists a $\Gal(K_{E, 13}/\Q)$-equivariant surjective homomorphism 
 \[
 \left(\cl(K_{E, 13})\tens{\Z}{\F_{13}}\right)^{\mathrm{ss}} \twoheadrightarrow E[13]
 \]
 by Corollary \ref{cor1}. In particular, we obtain
 \[
 r_{E[13]} \geqslant 1.
 \]
\end{example}

\begin{rem}\label{PS}
As we noted in the last of Subsection \ref{back}, Prasad and Shekhar proved that 
\[
\dim_{\F_p}(\Sel(\Q, E[p]))\geqslant 2 \Rightarrow  r_{E[p]} \neq 0
\]
under certain conditions in \cite[Theorem 3.1]{PS}. Due to LMDFB, the Tate-Shafarevich group $\Sha(E/\Q)$ of $E$ in Example \ref{ex} is trivial. Since $r(E)=1$ and $E(\Q)[13]=0$ by Mazur's theorem \cite[Chapter VIII, Theorem 7.5]{Sil}, we have $\dim_{\F_{13}}\Sel(\Q, E[13])=1$. Thus we cannot deduce that $r_{E[13]} \neq 0$ in this example with the result of Prasad and Shekhar. 
\end{rem}

\begin{example}\label{sugoi}
Let $p=7$ and $E$ be the elliptic curve defined by the equation 
\[
y^2=x^3-2401x+1.
\]
We apply Corollary \ref{cor1} to this example. This elliptic curve is a specialization at $m=49,\ n=1$ of the family of elliptic curves with two parameters $m, n$ treated in \cite{FN}. We can see that $E$ has three rational points 
\[
P:=(0, 1),\ Q:=(-49, 1)\ \text{and}\ R:=(-1, 49).
\]
 They are linearly independent and can be extended to a basis of $E(\Q)$ by \cite[Theorem 1.1 (2)]{FN}. In particular, they are not in $7E(\Q)$.

By a simple computation, we know that $j(E)=2^8\cdot 3^3\cdot 7^{12}/1069\cdot 51791533$, and (Tam) in Proposition \ref{mainprop} is satisfied by Remark \ref{tamj}. We can check that $E[7]$ is irreducible as a $\Gal(K_{E, 7}/\Q)$-module by MAGMA.

With the help of MAGMA, we obtain 
\[
3P=\left(\frac{-2^3\cdot 79 \cdot 199 \cdot 367 \cdot 2399}{{\bf 7}^{16}}, \frac{37\cdot 4691\cdot 19523423\cdot 169609859}{{\bf 7}^{24}}\right),
\]
\[
3Q=\left(\frac{5^2\cdot 13\cdot 53\cdot 181 \cdot 1777\cdot 73483}{2^2\cdot {\bf 7}^4 \cdot 67^2\cdot 439}, \frac{29\cdot 31\cdot 6151 \cdot 12992635846499}{2^3\cdot {\bf 7}^6\cdot 67^3\cdot 439^3} \right),
\]
\[
2R=\left(\frac{3^2\cdot 139\cdot 1153}{{\bf 7}^4}, \frac{5\cdot 345311039}{{\bf 7}^6}  \right).
\]
Then $3P, 3Q, 2R \in E(\Q)_{\ur, 7}\setminus 7E(\Q)$ by Theorem \ref{main2}. Since the three points $3P, 3Q, 2R$ give a basis of $E(\Q)\tens{\Z}{\Z_7}$ over $\Z_7$, their image in $E(\Q)/7E(\Q)$ give a basis of $E(\Q)/7E(\Q)$ over $\F_7$. In particular, they are linearly independent over $\F_7$ in $E(\Q)/7E(\Q)$. Thus we obtain
\[
r_{\ur, 7}(E) \geqslant 3.
\]
By Corollary \ref{cor1}, there exists a $\Gal(K_{E,7}/\Q)$-equivariant surjective homomorphsim
\begin{align}
\left(\cl(K_{E, 7})\tens{\Z}{\F_7}\right)^{\mathrm{ss}} \twoheadrightarrow E[7]^{\oplus 3}.
\end{align}
In particular, we obtain
\[
r_{E[7]} \geqslant 3.
\]
\end{example}
\begin{rem}
 In this example, \cite[Theorem 3.1]{PS} by Prasad and Shekhar gives the non-vanishing $r_{E[7]} \neq 0$ since we have $\dim_{\F_p}(\Sel(\Q, E[p])) \geqslant r(E) \geqslant 3$. Thus Corollary \ref{cor1} refines their result in this case.
\end{rem}

\end{document}